\documentclass[10pt,fleqn]{article}

\usepackage{graphicx}
\usepackage{amsmath}
\usepackage{amssymb}
\usepackage{amsthm}

\theoremstyle{plain}
\newtheorem*{thm}{Theorem}

\newtheorem*{cor}{Corollary}

\date{November 4, 2003}

\begin{document}

\title{Achievable ranks of intersections \\ of finitely generated free groups}
\author{Richard P. Kent IV}
\maketitle

\begin{abstract} We answer a question due to A. Myasnikov by proving that all expected ranks occur as the ranks of intersections of finitely generated subgroups of free groups.

\end{abstract}

\noindent Let $F$ be a free group.  Let $H$ and $K$ be nontrivial finitely generated subgroups of $F$.  It is a theorem of Howson \cite{how} that $H \cap K$ has finite rank.  H. Neumann proved in \cite{hanna} that $\mathrm{rank}(H \cap K) -1 \leq 2 (\mathrm{rank}(H)-1)(\mathrm{rank}(K)-1)$ and asked whether or not $\mathrm{rank}(H \cap K) -1 \leq (\mathrm{rank}(H)-1)(\mathrm{rank}(K)-1)$.

 A. Miasnikov has asked which values between $1$ and $(m-1)(n-1)$ can be achieved as $\mathrm{rank}(H\cap K) - 1$ for subgroups $H$ and $K$ of ranks $m$ and $n$---this is problem AUX1 of \cite{open}.  We prove that all such numbers occur by proving the following 

\begin{thm} Let $F(a,b)$ be a free group of rank two.  Let 
\begin{align*}
H_{k,\ell}^m =  \langle a,\  & bab^{-1}, \ldots,\  b^k a b^{-k},\  b^{k+1} a^{n-\ell} b^{-(k+1)}, &
\\ & b^{k+2} a^n b^{-(k+2)},\  b^{k+3} a^n b^{-(k+3)},  \ldots ,\ b^{m-1} a^n b^{1-m}  \rangle
\end{align*} 
and let $K = \langle b,\ aba^{-1}, \ldots ,\ a^{n-1}ba^{1-n}  \rangle$, where $0\leq k \leq m-2$ and $0 \leq \ell \leq n-1$.  Then the rank of $H_{k,\ell}^m \cap K$ is $k(n-1) + \ell$.
\end{thm}

\begin{cor} Let $F$ be a free group and let $m,n \geq 2$ be  natural numbers. Let $N$ be a natural number such that $1 \leq N - 1 \leq (m-1)(n-1)$.  Then there exist subgroups $H,K \leq F$, of ranks $m$ and $n$, such that the rank of $H \cap K$ is $N$.
\end{cor}
\begin{proof}[Proof of the corollary] The theorem produces the desired subgroups for all $N$ with $N-1 \leq (m-1)(n-1) - 1$ after passing to a rank two subgroup of $F$.  For $N -1 = (m-1)(n-1)$, simply let $H = \langle a,\ bab^{-1}, \ldots,\ b^{m-2}ab^{2-m},\ b^{m-1} \rangle$ and let $K = \langle b,\ aba^{-1}, \ldots,\ a^{n-2}ba^{2-n},\ a^{n-1} \rangle$.
\end{proof}

\begin{proof}[Proof of the theorem]  Let $X$ be a wedge of two circles and base $\pi_1(X)$ at the wedge point.  We identify $\pi_1(X)$ with $F=F(a,b)$ by calling the homotopy class of one oriented circle $a$ and the other $b$. Given a finitely generated subgroup of $F$, there is a covering space $\widetilde X$ corresponding to this subgroup.  Moreover, there is a compact subgraph of $\widetilde X$ that carries the given subgroup.  Given two subgroups and their associated finite graphs, one may construct the graph associated to their intersection.  These procedures are laid out carefully in \cite{stallings} and we assume that the reader is familiar with that paper.

 In the figures, the graph associated to $H$ appears at the top, that of $K$ to the right, and that of $H \cap K$ in the center. Edges labelled with two arrowheads represent $a$, those with one arrowhead represent $b$.  Our basepoint in the graph associated to $H \cap K$ is always the vertex in the upperlefthand corner. 

 For the moment, fix $k=m-2$. In Figure \ref{first}, $\ell=n-1$ and the rank of $H_{m-2,n-1}^m \cap K$ is visibly $(m-1)(n-1)$. Decreasing $\ell$ by one alters the intersection graph as depicted in Figure \ref{second} and the rank of $H_{m-2,n-2}^m \cap K$ is $(m-1)(n-1) - 1$. Figure \ref{third} shows the case when $\ell = n-3$ and the rank of the intersection is $(m-1)(n-1) - 2$.  When $\ell = n-j$, the rank of $H_{m-2,n-j}^m \cap K$ is $(m-1)(n-1) - (j-1)$. 

\begin{figure}
   \begin{center}
    \includegraphics{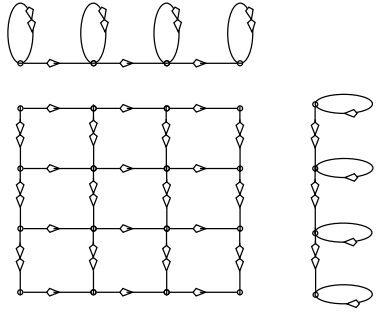}
   \end{center}
    \caption{$H$, $K$, and $H \cap K$ when $k=m-2$, $\ell = n-1$}
    \label{first}
    \end{figure}

\begin{figure}
   \begin{center}
    \includegraphics{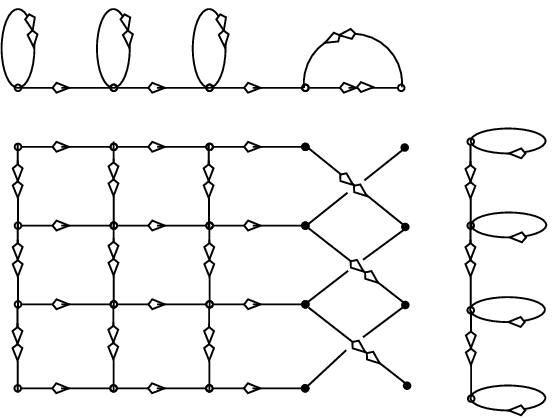}
   \end{center}
    \caption{$H$, $K$, and $H \cap K$ when $k=m-2$, $\ell = n-2$}
    \label{second}
    \end{figure}

\begin{figure}
   \begin{center}
    \includegraphics{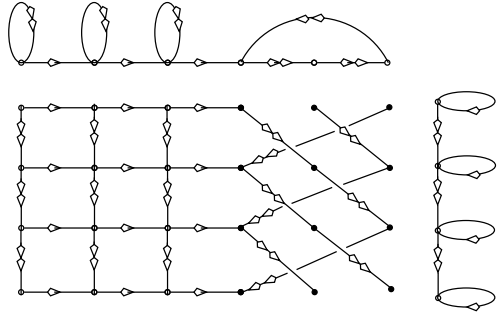}
   \end{center}
    \caption{$H$, $K$, and $H \cap K$ when $k=m-2$, $\ell = n-3$}
    \label{third}
    \end{figure}

 Figure \ref{fourth} depicts the case $\ell = 0$. Note that the graph associated to $H_{m-2,0}^m \cap K$ is the graph associated to $H_{m-3, n-1}^{m-1} \cap K$ to which a collection of trees have been attached at their roots, the graph associated to $H_{m-3,n-2}^m \cap K$ is the graph associated to $H_{m-3,n-2}^{m-1}\cap K$ to which trees have been so attached, and so on. Since attaching trees in this way leaves the rank intact, we arrive at the theorem by induction on $m$.
\begin{figure}
   \begin{center}
    \includegraphics{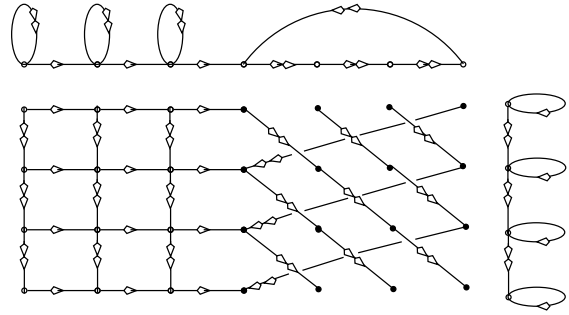}
   \end{center}
    \caption{$H$, $K$, and $H \cap K$ when $k=m-2$, $\ell = 0$}
    \label{fourth}
    \end{figure}
\end{proof}

\section*{Acknowledgement}

This work supported in part by a University of Texas Continuing Fellowship.

\noindent \textit{Department of Mathematics, University of Texas, Austin, TX 78712} 
\newline \noindent  \texttt{rkent@math.utexas.edu}

\end{document}